\documentclass[]{imsart}

\usepackage[latin1]{inputenc}

\usepackage{amsthm,amsmath}
\RequirePackage[colorlinks,citecolor=blue,urlcolor=blue]{hyperref}

\usepackage{hypernat}

\usepackage{amsfonts}%
\usepackage{amssymb}%
\usepackage{graphicx}
\usepackage{mathtools}
\usepackage{subcaption}
\DeclarePairedDelimiter{\ceil}{\lceil}{\rceil}

\startlocaldefs

\theoremstyle{plain}

\numberwithin{equation}{section}
\endlocaldefs

\begin{document}
\begin{frontmatter}

\title{ The relationship between stopping time and number of odd terms in Collatz sequences }
\runtitle{ Stopping time and odd terms in Collatz sequences }


\author{\fnms{Rafael} \snm{Ruggiero} \ead[label=e1]{ruggiero3n1@gmail.com}}
\address{S\~ao Paulo, Brazil\\
\printead{e1}}



\runauthor{R. Ruggiero}

\begin{abstract}
The Collatz sequence for a given natural number $N$ is generated by repeatedly applying the map $N$ $\rightarrow$ $3N+1$ if $N$ is odd and  $N$ $\rightarrow$ $N/2$ if $N$ is even. One elusive open problem in Mathematics is whether all such sequences end in 1 (Collatz conjecture), the alternative being the possibility of cycles or of unbounded sequences. In this paper, we present a formula relating the stopping time and the number of odd terms in a Collatz sequence, obtained numerically and tested for all numbers up to $10^7$ and for random numbers up to $2^{128.000}$. This result is presented as a conjecture, and with the hope that it could be useful for constructing a proof of the Collatz conjecture.
\end{abstract}





\end{frontmatter}

\section{Introduction}

The Collatz conjecture states that the sequence resulting from repeatedly applying the map $N$ $\rightarrow$ $3N+1$ if $N$ is odd and $N$ $\rightarrow$ $N/2$ if $N$ is even, starting from any natural number, inevitably reaches $1$. The alternative would be the existence of some special $N$ for which the sequence diverges or features an indefinite cycle. For a review on general properties of Collatz sequences, see \cite{lagarias}.

Despite the simplicity in the enunciation of the conjecture, it turns out to be an incredibly challenging problem. There is no simple formula relating $N$ to the number of iterations needed to reach $1$, and the behavior of this number (called the \emph{stopping time}) as a function of $N$ is markedly unpredictable, even passing several standard tests of randomness \cite{xu}.

Our goal in this paper is to propose an explicit formula for the stopping time $S$ as a function of the number of odd terms in the sequence and of $N$, a result which, if correct, allows one to constraint which values $S$ may assume.

\section{The formula}

Let $S$ be the stopping time for the Collatz map starting from a given natural number $N$, defined as the number of iterations needed for the sequence to reach $1$, and let $\alpha$ be the number of odd terms in the resulting sequence, excluding $1$.

\vspace{2.8cm}

\noindent \textbf{Conjecture:} For all $N$ for which $S$ is defined, it is true that
\begin{equation} \label{eq:main}
S = \ceil{\log_2\left(6^{\alpha} N\right)},
\end{equation}
where $\ceil{x}$ is the ceiling function, defined as the smallest integer greater than or equal to $x$. 

\begin{figure}[b!]
  \includegraphics[width=0.73\linewidth]{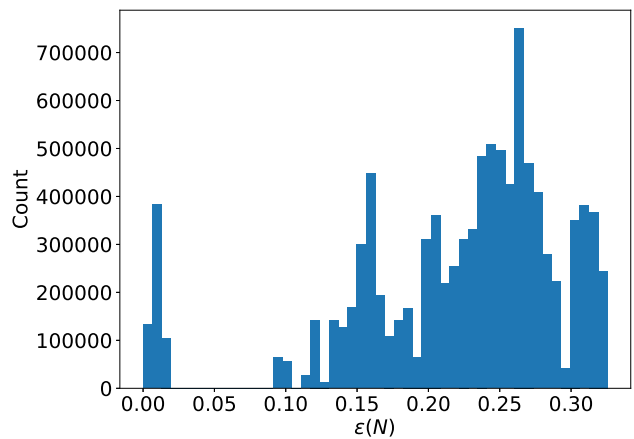}
  \caption{Distribution of the residue $\epsilon(N)$ for $N < 10^7$.}
  \label{fig:dist}
\end{figure}

The use of the ceiling function is necessary because the residue $\epsilon(N) \equiv S - \log_2\left(6^{\alpha} N\right)$ is in general not zero. In fact, it is a seemingly random real number in the interval $[0, 0.326)$. The distribution of $\epsilon(N)$ for all $N < 10^7$ is shown in Fig.~\ref{fig:dist}.

Our numerical tests have verified the validity of this relation for all natural numbers smaller than $10^7$, and also for thousands of random numbers between $1$ and $2^{128.000}$, strongly suggesting its validity for all natural numbers. A proof is left as an open problem.

\section{Implications}

In principle, we do not know how many odd numbers a Collatz sequence will have, so $\alpha$ is an unknown prior to having calculated the entire sequence. But if Eq.~\ref{eq:main} is valid, some values which $S$ \emph{cannot} assume can be calculated without calculating any iteration of the Collatz map.

First, since $\alpha \geq 0$, then $S \geq \log_2{N}$. This is a trivial result: the fastest a number could get to $1$ under the Collatz map is if the terms in its Collatz sequence never increase in value. In practice, this happens for powers of $2$.

A less trivial result is that, since $\alpha$ must be a natural number, not necessarily $S$ can assume any natural value $\geq \log_2{N}$. For instance, consider $N = 65$. Its possible stopping times, setting $\alpha$ as all possible natural numbers (including zero), are 7, 9, 12, 14, 17, 19, 22, 25, 27, 30, 32, 35, \textellipsis. The true value, which is $27$, is present in this list, but many numbers both smaller and larger than $27$ are not, say 15, 26 or 34. Let's call such numbers \emph{prohibited} values for $S(N)$.

The list of prohibited values $P$ for $S$ given a particular number $N$ can be propagated to infinite other numbers, as long as the Collatz conjecture is true. If that is the case, then $S(N/2) = S(N)-1$ if $N$ is even and $S(3N+1) = S(N)-1$ if $N$ is odd. Moreover, $S(2N) = S(N)+1$ for all $N$, and $S((N-1)/3) = S(N)+1$ if $(N-1)/3$ is an odd natural number. Using these relations, either $P+1$ or $P-1$ can also be shown to be prohibited values for these up to $3$ other numbers. This procedure can be carried out recursively and starting from different $N$ and associated $P$, resulting in a sieve algorithm for ruling out stopping times for natural numbers.

\vspace{0.31cm}
\section{Sequences with constant $\alpha$}

\begin{center}
\begin{figure}
\makebox[\linewidth][c]{
  \includegraphics[width=1.4\linewidth]{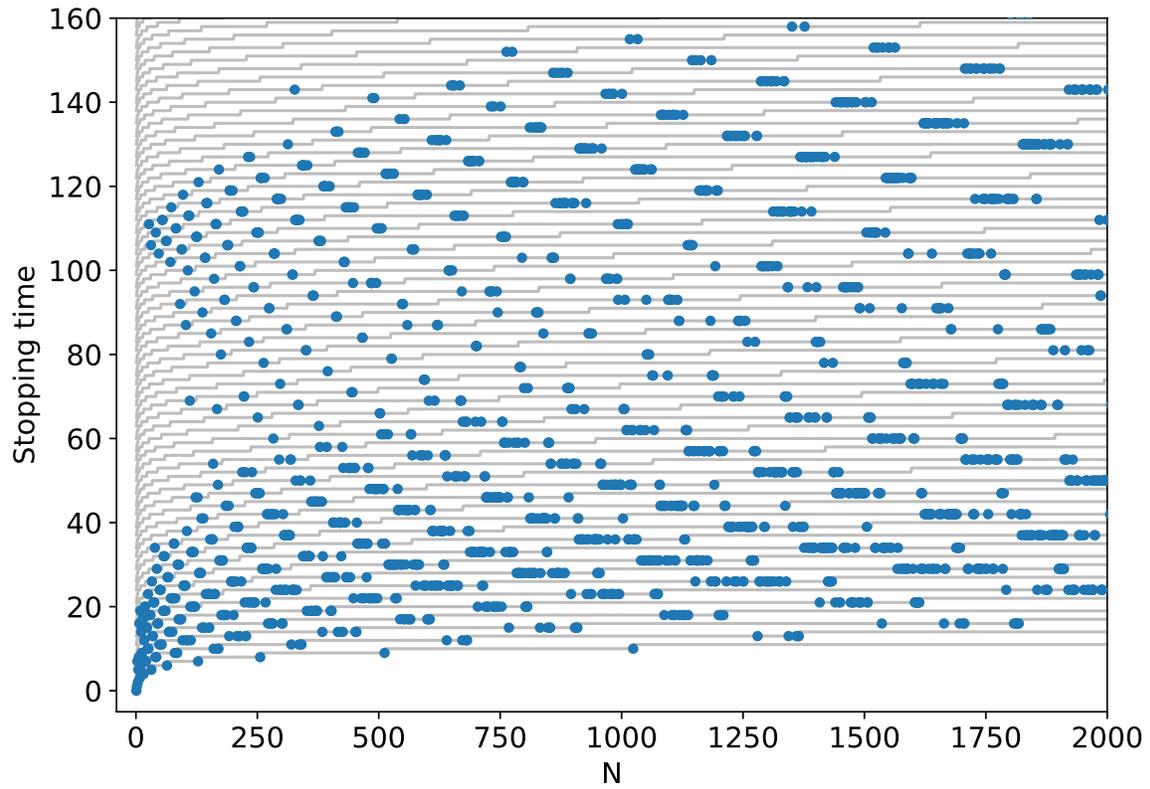}}
  \caption{Stopping time $S$ as a function of $N$ for $N < 2000$. The gray lines show Eq. \ref{eq:main} calculated for different natural values of $\alpha$ (including 0). The first line from bottom to top corresponds to $\alpha = 0$, the next one to $\alpha = 1$, and so on.}
  \label{fig:sn}
\end{figure}
\end{center}

\begin{figure}
\makebox[\linewidth][c]{
  \includegraphics[width=1.1\linewidth]{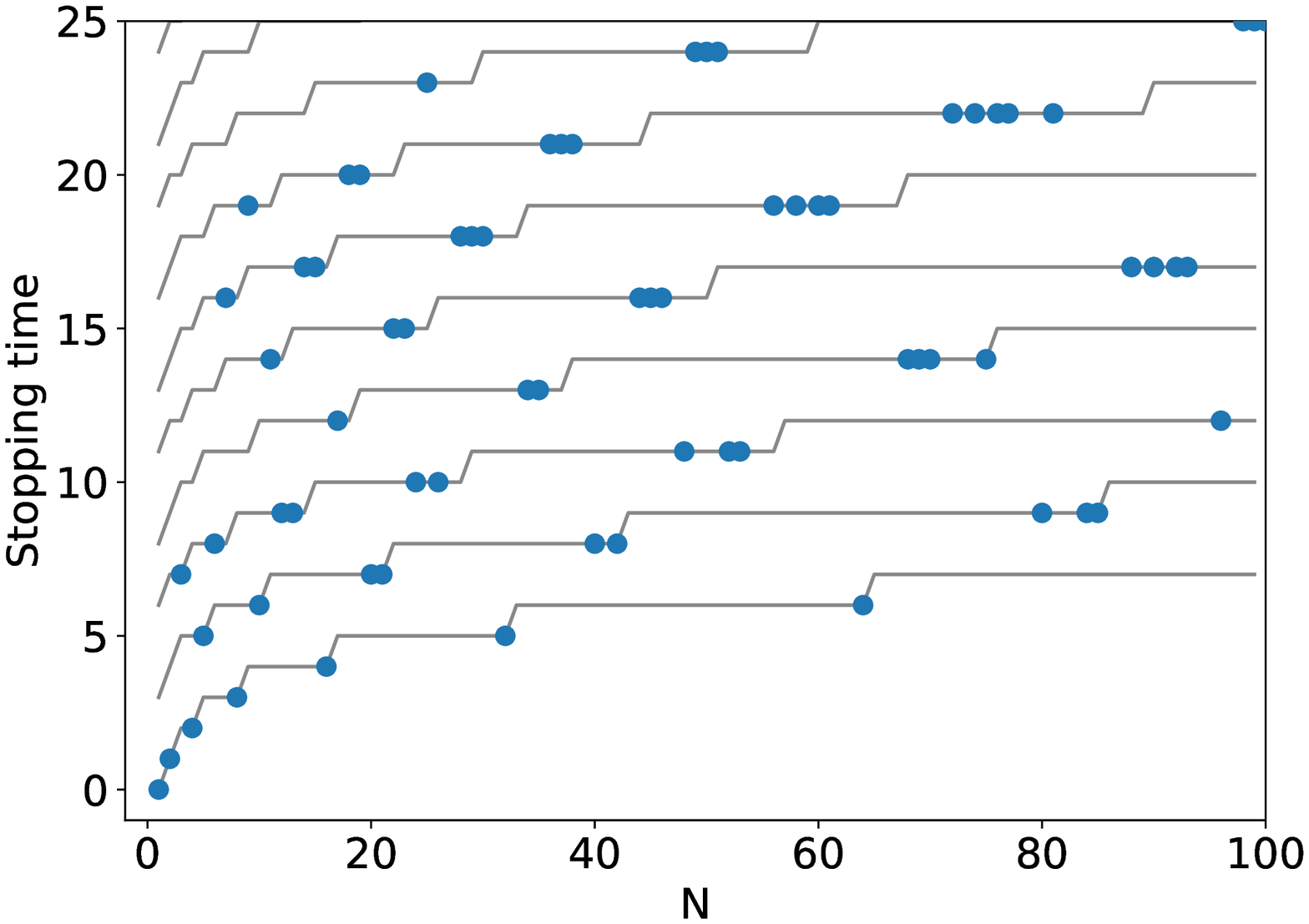}}
\makebox[\linewidth][c]{
  \includegraphics[width=1.1\linewidth]{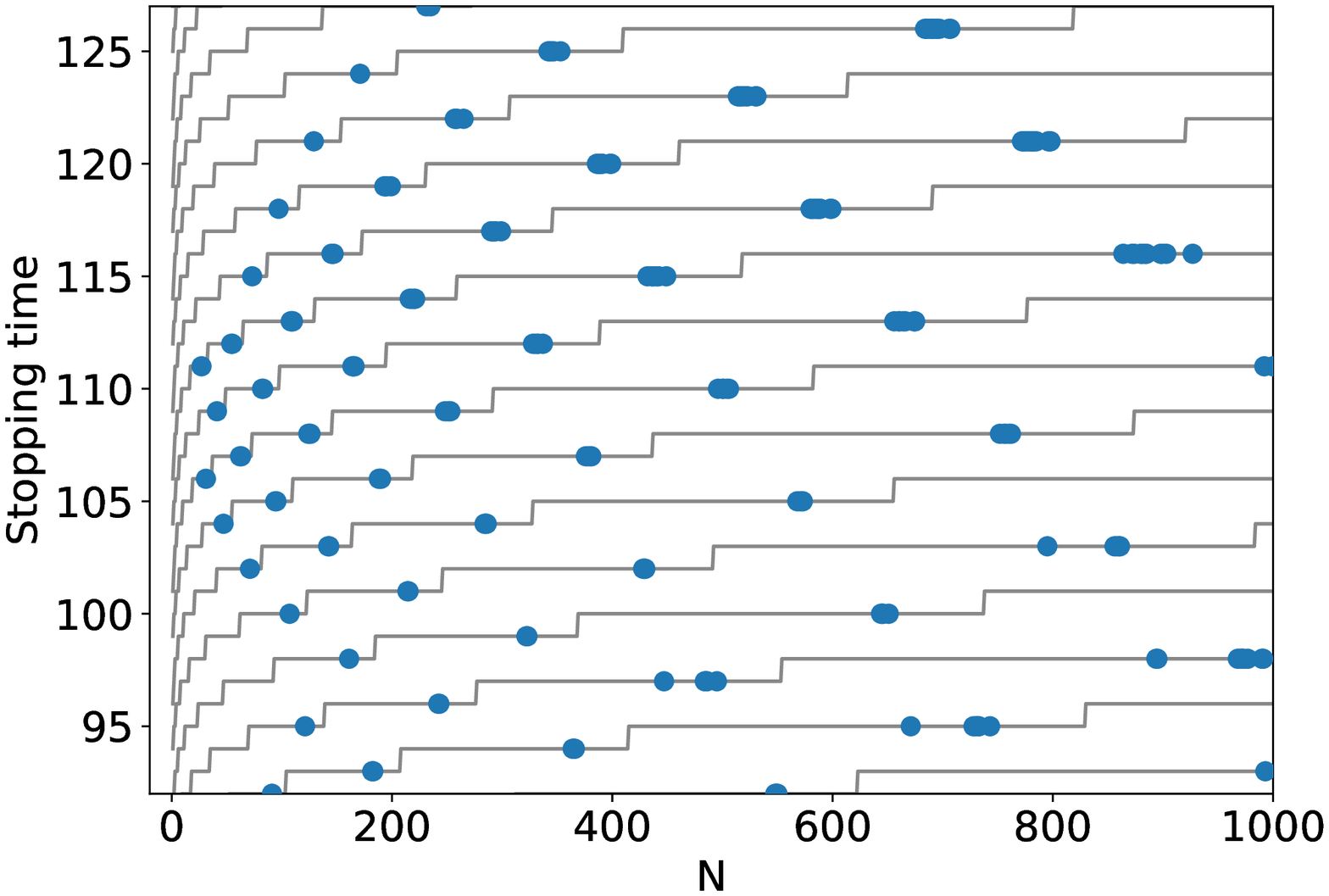}}
  \caption{Two zoomed-in rectangles in the area of Fig.~\ref{fig:sn}. The first one shows the points near the origin, and the second one shows the dense region on the upper left of that plot.}
  \label{fig:z}
\end{figure}

One important feature of Eq.~\ref{eq:main} is that it connects the stopping times of sequences of natural numbers with constant $\alpha$ in an $S(N)$ plot, such as the one shown in Fig.~\ref{fig:sn}, in which stopping times and curves with constant $\alpha$ are overlaid. The first such sequence, which is the lower frontier in that figure, is the one for $\alpha = 0$. This is the sequence of all powers of $2$, which have no odd numbers in their Collatz sequences. Sequences with higher $\alpha$ appear increasingly higher, and no two curves with different $\alpha$ ever intersect. We have also included two zoomed-in versions of Fig.~\ref{fig:sn} in Fig.~\ref{fig:z}, and the first elements of some sequences of numbers with constant $\alpha$ are shown in Table~\ref{tab:tab}.

\begin{table} [h!]
\caption{\normalsize The first elements of the sequences of numbers with constant $\alpha$, for $\alpha < 20$.}
\vspace{0.3cm}
\normalsize
\makebox[\linewidth][c]{
\begin{tabular}{ l l }
 $\alpha$ = 0: & 1, 2, 4, 8, 16, 32, 64, 128, 256, 512, 1024, 2048, 4096, 8192, 16384, 32768, 65536 \\
$\alpha$ = 1: & 5, 10, 20, 21, 40, 42, 80, 84, 85, 160, 168, 170, 320, 336, 340, 341, 640 \\
$\alpha$ = 2: & 3, 6, 12, 13, 24, 26, 48, 52, 53, 96, 104, 106, 113, 192, 208, 212, 213 \\
$\alpha$ = 3: & 17, 34, 35, 68, 69, 70, 75, 136, 138, 140, 141, 150, 151, 272, 276, 277, 280 \\
$\alpha$ = 4: & 11, 22, 23, 44, 45, 46, 88, 90, 92, 93, 176, 180, 181, 184, 186, 201, 352 \\
$\alpha$ = 5: & 7, 14, 15, 28, 29, 30, 56, 58, 60, 61, 112, 116, 117, 120, 122, 224, 232 \\
$\alpha$ = 6: & 9, 18, 19, 36, 37, 38, 72, 74, 76, 77, 81, 144, 148, 149, 152, 154, 162 \\
$\alpha$ = 7: & 25, 49, 50, 51, 98, 99, 100, 101, 102, 196, 197, 198, 200, 202, 204, 205, 217 \\
$\alpha$ = 8: & 33, 65, 66, 67, 130, 131, 132, 133, 134, 260, 261, 262, 264, 266, 268, 269, 273 \\
$\alpha$ = 9: & 43, 86, 87, 89, 172, 173, 174, 177, 178, 179, 344, 346, 348, 349, 354, 355, 356 \\
$\alpha$ = 10: & 57, 59, 114, 115, 118, 119, 228, 229, 230, 236, 237, 238, 456, 458, 460, 461, 465 \\
$\alpha$ = 11: & 39, 78, 79, 153, 156, 157, 158, 305, 306, 307, 312, 314, 315, 316, 317, 610, 611 \\
$\alpha$ = 12: & 105, 203, 209, 210, 211, 406, 407, 409, 418, 419, 420, 421, 422, 431, 455, 812, 813 \\
$\alpha$ = 13: & 135, 139, 270, 271, 278, 279, 281, 287, 303, 540, 541, 542, 545, 551, 556, 557, 558 \\
$\alpha$ = 14: & 185, 187, 191, 361, 363, 367, 370, 371, 374, 375, 382, 383, 721, 722, 723, 726, 727 \\
$\alpha$ = 15: & 123, 127, 246, 247, 249, 254, 255, 481, 489, 492, 493, 494, 498, 499, 508, 509, 510 \\
$\alpha$ = 16: & 169, 329, 338, 339, 359, 641, 657, 658, 659, 665, 676, 677, 678, 718, 719, 1281, 1282 \\
$\alpha$ = 17: & 219, 225, 239, 427, 438, 439, 443, 450, 451, 478, 479, 854, 855, 876, 877, 878, 886 \\
$\alpha$ = 18: & 159, 295, 318, 319, 569, 585, 590, 591, 601, 636, 637, 638, 1138, 1139, 1147, 1151, 1159 \\
$\alpha$ = 19: & 379, 393, 425, 758, 759, 767, 779, 786, 787, 801, 849, 850, 851, 1516, 1517, 1518, 1529 \\
\end{tabular}
}\label{tab:tab}
\end{table}

\vspace{-0.37cm}
\section{Collatz sequences in an $S(N)$ plot}

\begin{figure}
  \includegraphics[width=0.75\linewidth]{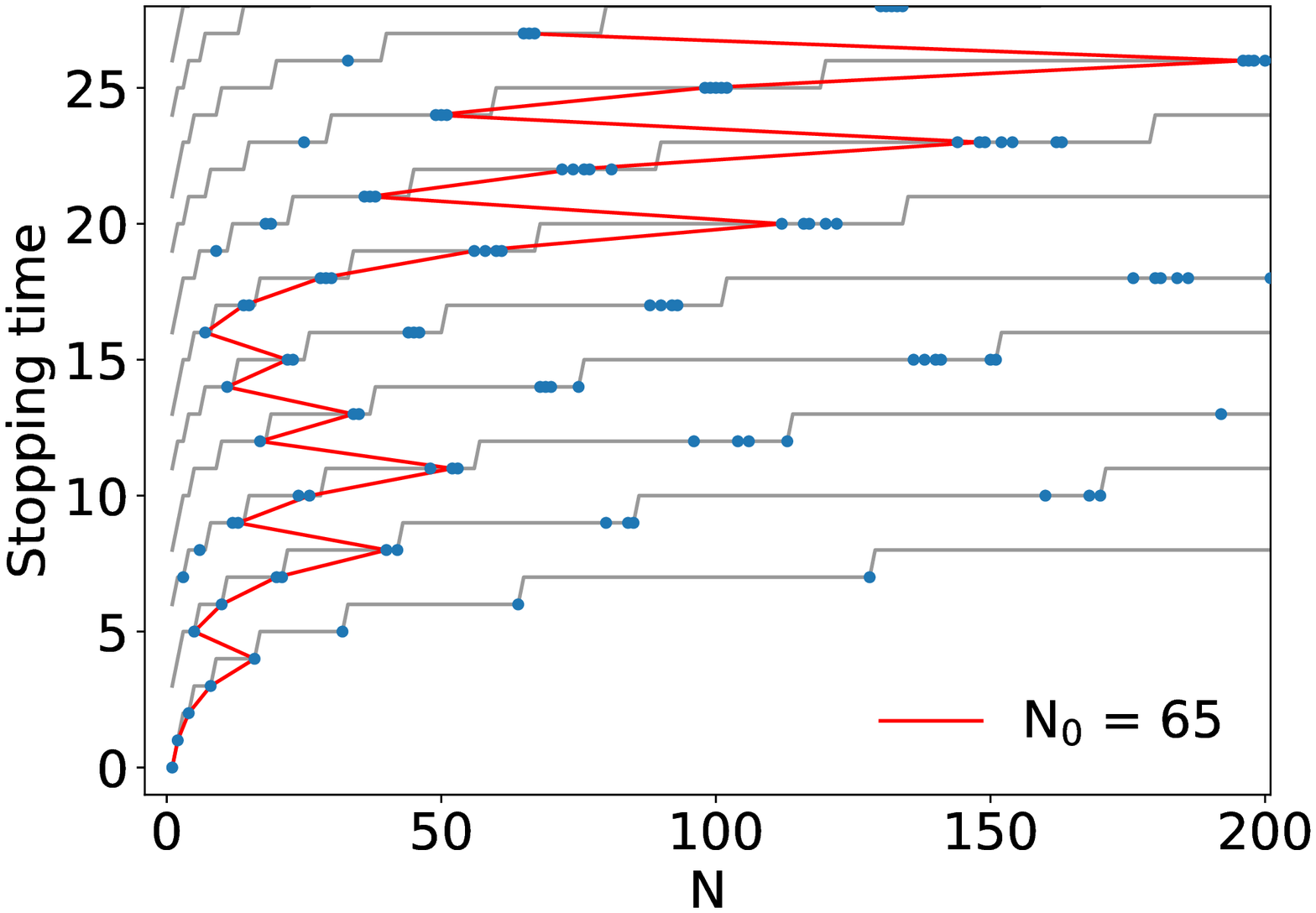}
  \includegraphics[width=0.75\linewidth]{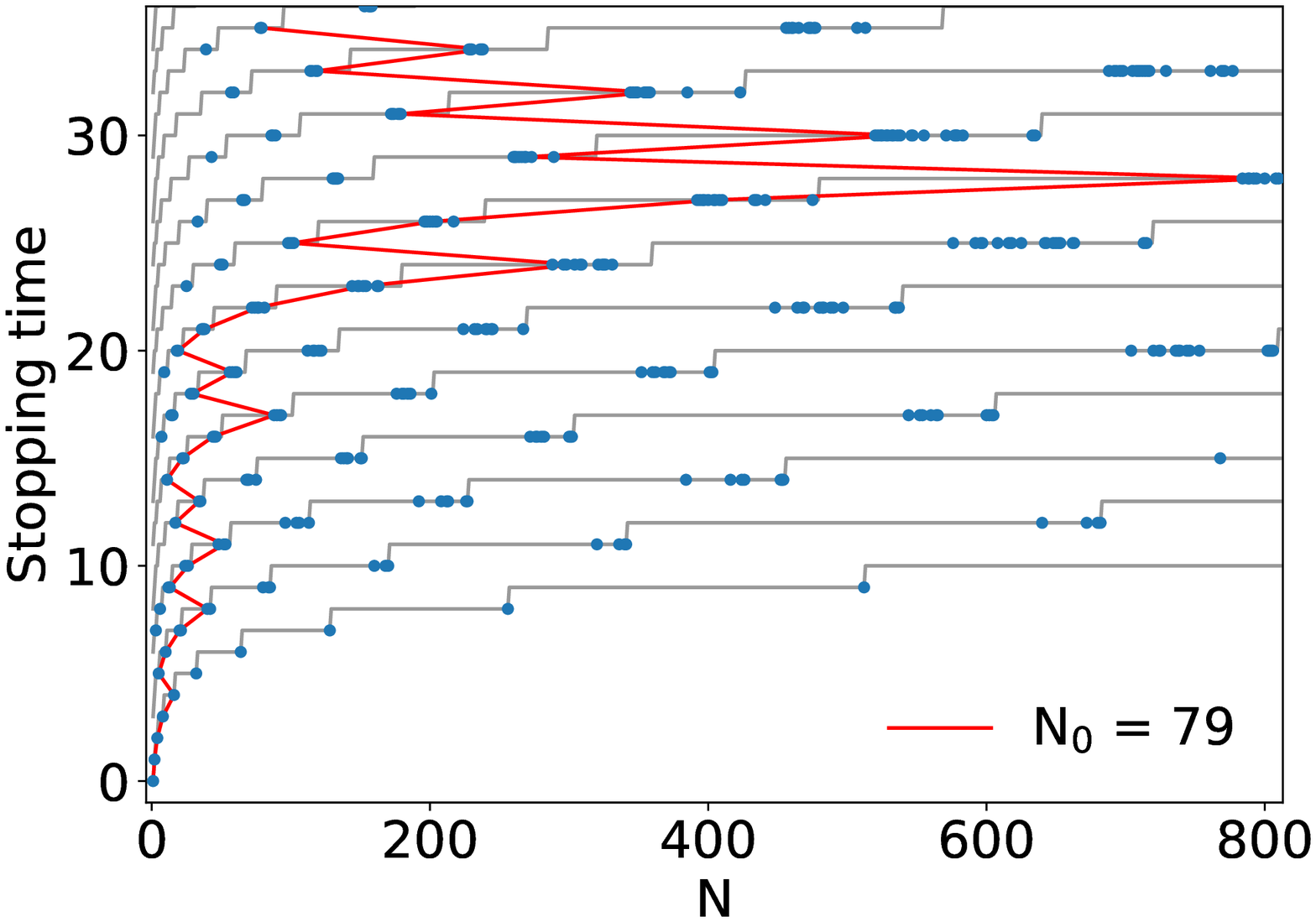}
  \includegraphics[width=0.75\linewidth]{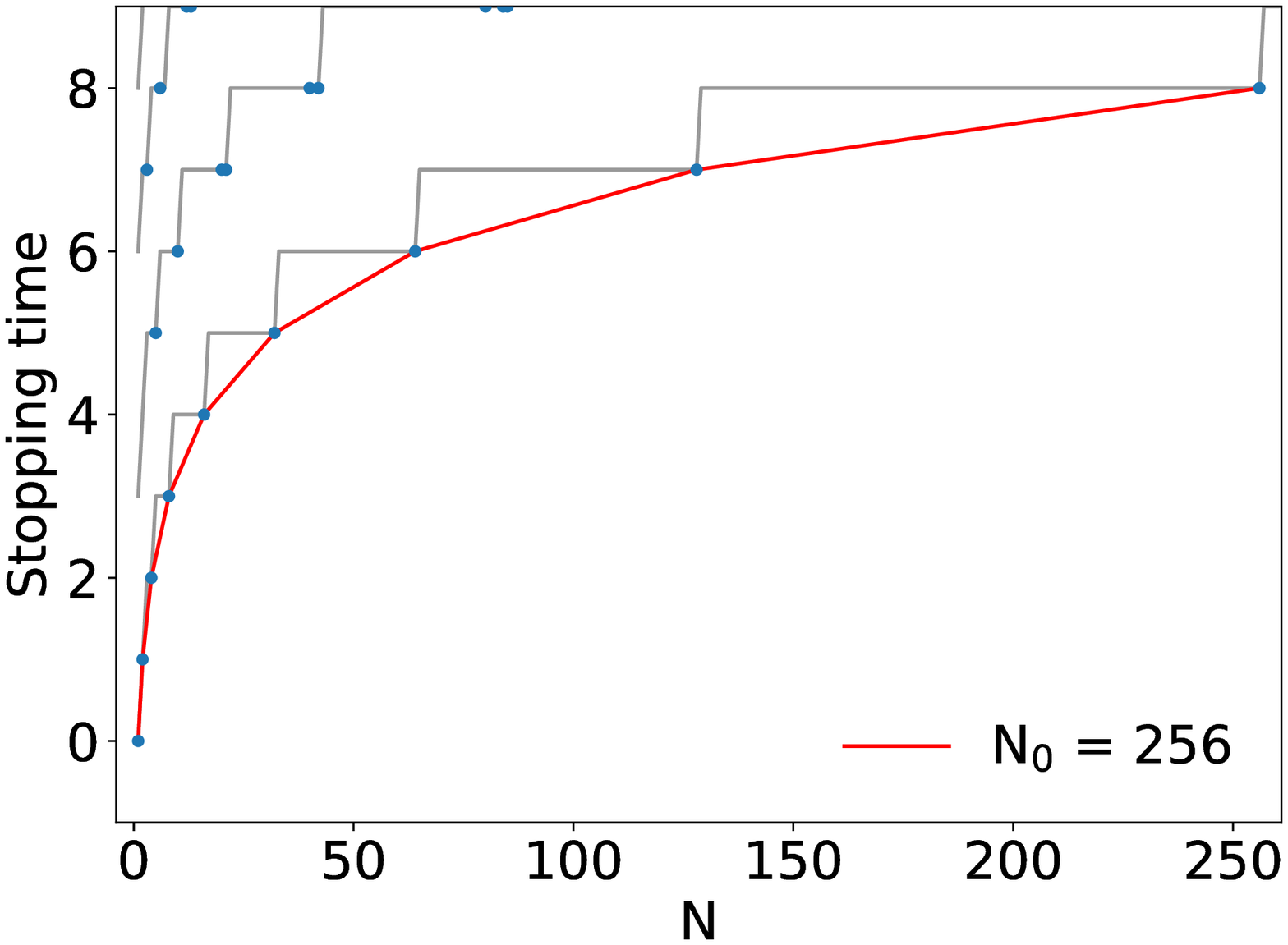}
  \caption{Examples of Collatz sequences in $S(N)$ plots, starting from different initial values. The red lines begin on top and move lower towards $N = 1$ at each iteration. The number either moves to the left along a curve with constant $\alpha$ or to the right to a point in a curve with smaller $\alpha$.}
  \label{fig:examples}
\end{figure}

It is illustrative to show what a Collatz sequence looks like in an $S(N)$ plot. As the sequence goes on, $S(N)$ decreases by 1 at each step, and one of two possibilities takes place: either the point moves along a curve with constant $\alpha$ to the left, as its factors of two are removed, or it moves right to a curve with lower $\alpha$ underneath the previous one. The process is repeated until the point becomes a power of $2$ and lands at the curve with $\alpha = 0$, from where it monotonically decreases towards $1$.
Three examples of this process are shown in Fig.~\ref{fig:examples}.

\end{document}